\documentclass{amsart}
\usepackage{amscd}
\DeclareMathOperator*{\colim}{colim}
\begin{document}
\newtheorem{prop}{Proposition}
\newtheorem{thm}[prop]{Theorem}
\newtheorem{ex}[prop]{Example}
\newtheorem{lemma}[prop]{Lemma}
\newtheorem{defn}[prop]{Definition}
\newtheorem{cond}[prop]{Condition}
\newtheorem{cor}[prop]{Corollary}
\newcommand{\into}{\hookrightarrow}

\newcommand{\comment}[1]{}
\newtheorem{remark}[prop]{Remark}

\newcommand{\bde}{\begin{defn}}
\newcommand{\ede}{\end{defn}}

\newcommand{\bpr}{\begin{prop}}
\newcommand{\epr}{\end{prop}}
\newcommand{\bre}{\begin{remark}}
\newcommand{\ere}{\end{remark}}

\newcommand{\ar}[1]{\stackrel{#1}{\to}}

\newcommand{\CC}{\mathcal{C}}
\newcommand{\Z}{\mathbb{Z}}

\newcommand{\hh}{\mathbb{H}}
\newcommand{\R}{\mathbb{R}}
\newcommand{\N}{\mathbb{N}}
\newcommand{\Q}{\mathbb{Q}}
\newcommand{\io}{\iota}
\newcommand{\Om}{\Omega}
\newcommand{\Si}{\Sigma}
\newcommand{\La}{\Lambda}
\newcommand{\x}{\times}
\newcommand{\ot}{\otimes}
\newcommand{\G}{\mathcal{G}}
\newcommand{\sm}{\wedge}
\newcommand{\HH}{\mathbb{H}}
\newcommand{\ep}{\varepsilon}
\newcommand{\CP}{\mathbb{CP}}
\newcommand{\RP}{\mathbb{RP}}

\title{Generalized string topology operations}
\author{Kate Gruher}
\thanks{The first author was supported by a National Defense Science and Engineering Graduate Fellowship.}
\address{Department of Mathematics, Stanford University, Stanford, CA 94305}
\email{kagruher@stanford.edu}
\author{Paolo Salvatore}
\address{Dipartimento di matematica, Universit\`a di Roma ``Tor Vergata'', 00133 Roma, Italy}
\email{salvator@mat.uniroma2.it}
\date{}

\begin{abstract}
We show that the Chas-Sullivan loop product, a combination of the Pontrjagin product on the fiber and intersection product on the base,
 makes sense on the total space
homology of any fiberwise monoid $E$ over a closed oriented manifold $M$.
More generally the Thom spectrum $E^{-TM}$ is a ring spectrum.
Similarly a fiberwise module over $E$ defines a module over $E^{-TM}$.
Fiberwise monoids include adjoint bundles of principal bundles, and the construction is natural with respect to maps of principal bundles.  This naturality implies homotopy invariance of the
algebra structure on $H_*(LM)$ arising from the loop product. 
If $M=BG$ is the infinite dimensional classifying space of a compact Lie group,
then we get a well-defined pro-ring spectrum, which we define to be the string topology of $BG$. 
If $E$ has a fiberwise action of the little $n$-cubes operad then $E^{-TM}$ is
an $E_n$-ring spectrum. This gives homology operations combining Dyer-Lashof operations on the fiber and Steenrod operations
on the base.
We give several examples where the new operations
give homological insight,
borrowed from
knot theory, complex geometry, gauge theory, and homotopy theory.

\end{abstract}
\maketitle

\section{Introduction}

Let $M$ be a smooth, closed, oriented 
 manifold and let $LM$ be its loop space.
    Chas and Sullivan constructed in \cite{cs} new operations on the
   homology of $LM$.  This started a stream of ideas that now form a subject called ``string topology.''  The construction of the main operation, the loop product, mixes the intersection product on the homology of the manifold with the Pontrjagin product of its based loop space.
    We show that this construction can be extended in many directions.

First of all, the product can be constructed on the homology
 of the total space of any fiberwise monoid over a smooth, closed, oriented manifold.
   A fiberwise monoid is a fiber bundle such that each fiber is a monoid and the product varies continuously.  Cohen and Jones defined in \cite{cj} a ring spectrum structure on the Thom spectrum $LM^{-TM}$  which realizes the Chas-Sullivan product in homology, where $TM \to M$ is the tangent bundle of $M$ and
  the virtual bundle $-TM$ is pulled back over $LM$ via the evaluation $ev:LM \to M.$
 The Cohen-Jones construction can also be generalized to fiberwise monoids, yielding a ring spectrum.  Hence we describe a large class of ring spectra.  Notice that this generalized construction encompasses not only the Cohen-Jones ring spectrum structure on $LM^{-TM}$ but also the ring spectrum structure on $M^{-TM}$ arising from Atiyah duality.
 Analogous constructions work for fiberwise modules over fiberwise monoids.

Furthermore,  we will construct homology operations
on the total space homology of a bundle over $M$ with a fiberwise
action of the little $n$-cubes operad.
We shall relate such operations to Steenrod operations on the base
 and Dyer-Lashof operations on the fiber.
The extended Cohen-Jones construction produces
in this setting $E_n$-ring spectra
and modules over them.

Naturality properties for the extended Cohen-Jones construction for fiberwise monoids will allow us to define the notion of the string topology of the classifying space $BG$ of a compact lie group $G$,
 which is closely related to the loop space $LBG$.
  Since $BG$ is not a finite-dimensional manifold, this is not covered by the
  work of Chas-Sullivan or Cohen-Jones.  Indeed, we will see that the string topology
   structure does not appear on homology or on a ring spectrum,
   as in Chas-Sullivan and Cohen-Jones, but rather on pro-homology and a pro-ring spectrum, which are inverse systems of the structures arising from submanifolds of $BG$.  Throughout this construction
 we will utilize the fact that $LBG$ is homotopy equivalent
  to the adjoint bundle of the universal principal $G$ bundle, $EG$.

The naturality properties of the extended Cohen-Jones construction also yield a fairly elementary proof that if $f: M_1 \to M_2$ is a homotopy equivalence of smooth, closed, oriented manifolds, then $Lf_*: H_*(LM_1) \to H_*(LM_2)$ is a ring isomorphism.  This result was originally announced by Cohen, Klein and Sullivan \cite{cks}.  

We present a wide range of examples:
At the level of Thom spectra,
the space of smooth maps from a $n$-manifold $N$ to $M$
 gives a module spectrum
over the $E_n$-ring spectrum  associated to the space of smooth $n$-spheres into $M$. 
The space of knots in the 3-sphere gives an $E_2$-ring spectrum.
The space of rational curves on a homogeneous projective variety
 gives also an $E_2$-ring spectrum.
The classifying spaces of gauge groups give pro-ring and pro-module spectra.
Maps from Riemann surfaces of all genera into $M$ define a ring spectrum.
The ``string topology'' operations of these examples
allow new homological computations
\cite{SK,SS}.

The paper will be organized as follows:  
Section \ref{sec:fibmon} will be devoted to defining fiberwise monoids and modules.  We will give one of our main examples of fiberwise monoids - adjoint bundles of principal bundles - and define the product on the total space homology.  

In Section \ref{sec:ringspectra} we will generalize the Cohen-Jones construction to fiberwise monoids and modules, and see that the product from Section 2 defines a graded algebra structure on the total space homology.  In Section \ref{sec:nat} we will give some naturality properties for the ring spectra arising from fiberwise monoids and modules.  

Section \ref{sec:c_n} will give the construction of the homology operations for a bundle with fiberwise action of the little $n$-cubes operad.  Section \ref{sec:bg} will define the string topology of $BG$, and Section \ref{sec:hominv} will prove the homotopy invariance of the algebra structure on $H_*(LM).$  Section \ref{examples} will present numerous examples.

The first author would like to thank R. Cohen for his ideas, advice and guidance.  The second author is grateful to R. Cohen and S. Kallel
 for several fruitful discussions.

\section{Fiberwise Monoids and Modules}\label{sec:fibmon}

Fiberwise homotopy
 theory has been developed by James and Crabb \cite{crj}. 
  In this section we define fiberwise monoids, using a slight variation on their terminology, and show that a fiberwise monoid structure on a fiber bundle over a closed, oriented $d$-dimensional manifold induces a product of degree $-d$ on the homology of the total space of the bundle.

Given a locally trivial fiber bundle $E \stackrel{p}{\to} B$ with fiber $F$, let $E \times _B E$ denote the fiberwise product, 
$$E\times _B E = \{(x,y) \in E\times E : p(x)=p(y) \}.$$ 
Then $E\times _B E$ is clearly a fiber bundle over $B$ with fiber $F\times F.$ The projection $E\times _B E \to B$ will also be denoted $p$.  Given two fiber bundles $E_1\stackrel{p_1}{\to} B$ and $E_2 \stackrel{p_2}{\to} B$, a map $f: E_1 \to E_2$ is said to be fiberwise if $p_2 \circ f = p_1.$

\begin{defn}  A fiberwise monoid consists of a fiber bundle $E \stackrel{p}{\to} B$ along with a fiberwise map $m: E\times _B E \to E$ and a section $s:B \to E$ satisfying 
$$m(x, s(p(x)))=m(s(p(x)),x) = x$$
and 
$$m(m(x,y),z)=m(x,m(y,z))$$ 
for any $x,y,z \in E$. The map $m$ is called the multiplication map and $s$ is called the unit section.
\end{defn}

\emph{Remark.} For our constructions it is generally enough to have fiberwise homotopies (homotopies that are fiberwise maps over $B$ at each stage) instead of equalities in the definition above.  However, for simplicity we will assume equality holds and leave the generalizations to the fiberwise homotopy case to the reader.
\medskip

 A morphism $f:E_1 \to E_2$ of fiberwise monoids is a commutative diagram 
$$\begin{CD} 
E_1 @>f>> E_2 \\
@V{p_1}VV @V{p_2}VV \\
M_1 @>f>> M_2 
\end{CD}$$
that respects the multiplication and the unit section, i.e., $f\circ s_1 = s_2 \circ f$ and $f(m_1(x,y)) = m_2(f(x),f(y)).$

For a closed manifold $M$, the free loop bundle $LM \stackrel{ev}{\to} M$ has a fiberwise multiplication defined by $m(\alpha, \beta) = \alpha \ast \beta$ where
$$\alpha \ast \beta (t) = \begin{cases} \alpha(2t) & \quad 0\leq t \leq \frac{1}{2} \\ \beta(2t-1) & \quad \frac{1}{2}\leq t \leq 1. \end{cases}$$  
Now, $m$ is not strictly associative and thus does not make $LM$ into a fiberwise monoid.  However, $m$ is associative up to fiberwise homotopy.  Likewise the section $s:M\into LM$ that embeds $M$ as the constant loops is a unit up to fiberwise homotopy ($LM$ is what is called a fiberwise monoid-like space). As noted above, this is enough for our constructions.   Alternatively, replacing loops with ``Moore loops'' gives a bundle over $M$ that is an honest fiberwise monoid.

Adjoint bundles will provide our main example of fiberwise monoids.  Recall that a principal $G$ bundle $E \stackrel{p}{\to} B$ has a free, fiberwise right $G$ action.  We can form the adjoint bundle $Ad(E) \stackrel{\pi}{\to} B$ by 
$$Ad(E) = E\times _G G = E\times G/ \sim$$
where $(x,g)\sim (x,g)h = (xh, h^{-1}gh)$ for any $h$ in $G.$   Let $[x,g]$ denote the equivalence class of $(x,g)$; then $\pi([x,g]) = p(x).$  The adjoint bundle has fiber $G$ but in general is not a principal bundle.  If $G$ is abelian then the adjoint bundle is trivial, $Ad(E) \cong B\times G.$

The group structure on $G$ induces a fiberwise monoid structure on $Ad(E)$ as follows.  If $[x,g]$ and $[y,h]$ are in the same fiber of $Ad(E)$, then $x$ and $y$ are in the same fiber of $E$, so there exists a unique element $k \in G$ such that $y=xk.$  Then $[y,h]=[xk,h]=[x,khk^{-1}].$  Define $m([x,g],[y,h]) = [x,gkhk^{-1}].$ It is easy to check that $m$ is well-defined, continuous, fiberwise, and associative.  Notice that $[x,1] = [y,1]$ for any $x$ and $y$ in the same fiber of $E$.  Hence the unit section is given by $s(z) = [x,1]$ where $x$ is any element of $p^{-1}(z).$  Thus $m$ gives $Ad(E)$ a fiberwise monoid structure.  

In order to define the product on the homology of the total space of a fiberwise monoid, let us first recall the Pontrjagin-Thom construction.  For the usual Pontrjagin-Thom construction, let $e:P\into N$ be a smooth embedding of smooth manifolds.  Then there exists a tubular neighborhood $\eta_e$ of $e$, and $\eta_e$ is homeomorphic to the total space of the normal bundle $\nu_e.$  Then the Pontrjagin-Thom collapse map is 
$$\tau: N \to N/(N - \eta_e).$$
Identifying $N/(N - \eta_e)$ with the Thom space of the normal bundle $P^{\nu_e}$, we obtain a map $\tau: N \to P^{\nu_e}$.  

However, this Pontrjagin-Thom construction is sometimes not general enough, for instance when we want to work with infinite dimensional manifolds.  The following generalization is due to Klein and R. Cohen \cite{ralph}.  Let 
$$\begin{CD}
E @>>> X \\
@VqVV @VVpV \\
P @>e>> N 
\end{CD}$$
be a homotopy Cartesian square, where $e$ is a smooth embedding of closed manifolds.  Then there exists a Pontrjagin-Thom collapse map
$$X \to E^{q^*\nu_e}$$
well-defined up to homotopy.  Cohen's construction involves replacing $p$ by the homotopy equivalent path fibration $\tilde{X} \stackrel{\tilde{p}}{\to} N$, then showing that you can pick a tubular neighborhood $\eta_e$ of $e$ such that $\tilde{p}^{-1}(\eta_e)$ is homotopy equivalent to the total space of $q^*\nu_e.$  This collapse map is compatible with the collapse map coming from $e$ in the sense that 
$$\begin{CD}
X @>{\tau}>> E^{q^*\nu_e} \\
@VpVV @VVqV \\
N @>{\tau}>> P^{\nu_e}
\end{CD}$$
commutes.  Furthermore, in the case that $q: X \to N$ is a fiber bundle and $E = e^*X$, then there is no need to replace $X$ in the construction: $p^{-1}(\eta_e)$ is a tubular neighborhood of $E$ in $X$ and is homeomorphic to the total space of $q^*\nu_e$.    

Let $F\into E\stackrel{p}{\to} M$ be a fiberwise monoid with multiplication map $m$.  Assume $M$ is a closed, smooth, oriented manifold of dimension $d$.  Let $\Delta$ denote the diagonal map $\Delta: M \to M\times M$ and $\tilde{\Delta}$ the obvious codimension $d$ embedding of $E\times _M E$ into $E\times E$.  Then 
$$\begin{CD}
E\times _M E @>{\tilde{\Delta}}>> E\times E \\
@V{p}VV @VV{p\times p}V \\
M @>>\Delta > M\times M 
\end{CD}$$
commutes and is a pullback diagram of fiber bundles. Then as above there is a Pontrjagin-Thom collapse map $\tau: E \times E \to (E \times _M E)^{p^*\nu_{\Delta}}$.  Furthermore, the normal bundle of the diagonal embedding of a closed manifold is isomorphic to the tangent bundle, $\nu_\Delta \cong TM.$  Hence we obtain the collapse map 
$$\tau: E\times E \to (E\times_M E) ^{p^*TM}.$$
  For convenience we will often write the Thom space $(E \times _M E)^{p^*TM}$ as $(E\times _M E) ^{TM}$ since it is clear what map we are using to pull back the tangent bundle.

In homology, combining $\tau _*$ with $m_*$ and with the Thom isomorphism $$u_*:\tilde{H}_n((E\times_M E)^{TM})\stackrel{\cong}{\to} H_{n-d}(E\times _M E)$$ 
yields a product
\begin{multline*}
\mu_*: H_s(E)\otimes H_t(E) \stackrel{\times}{\to} H_{s+t}(E\times E) \stackrel{\tau_*}{\to} \\  \tilde{H}_{s+t}((E\times _M E)^{TM})
 \stackrel{u_*}{\to} \tilde{H}_{s+t-d}(E\times _M E) \stackrel{m_*}{\to} \tilde{H}_{s+t-d}(E).
\end{multline*}
In the case $E=LM$, this is the Chas-Sullivan product, as shown in \cite{cj}.  To simplify the grading, define $\mathbb{H}_*(E) = H_{*+d}(E)$, so that $\mu_*:\mathbb{H}_s(E)\otimes \mathbb{H}_t(E) \to \mathbb{H}_{s+t}(E).$  We will see in Section \ref{sec:ringspectra} that the product $\mu _*$ makes $\mathbb{H}_*(E)$ into a graded algebra, compatible with the algebra structure on $\mathbb{H}_*(M)$ arising from the intersection product.

Many of our constructions for fiberwise monoids easily generalize to fiberwise modules.
\begin{defn}
Let $E\stackrel{p}{\to} M$ be a fiberwise monoid and let $E' \to M$ be a fiber bundle along with a fiberwise map $m': E\times_M E' \to E'$ satisfying:
$$m'(x,m'(y,z)) = m'(m(x,y),z)$$
where $m$ is the multiplication map for $E \to M$, and 
$$m'(s(p(z)),z) = z$$
where $s$ is the unit section for $E\to M.$  Then $E'\to M$ is said to be a fiberwise module over $E\to M$.
\end{defn}

If $E' \to M$ is a fiberwise module over $E\to M$, where $M$ is a smooth, closed, oriented $d$-dimensional manifold, then we can construct a product $\mu ' _* :H_*(E) \otimes H_*(E') \to \tilde{H}_{*-d}(E')$ by composing the Pontrjagin-Thom collapse map $\tau$ and the Thom isomorphism with $m'_*$ analogous to the fiberwise monoid situation:
\begin{multline*}
\mu'_*: H_s(E) \otimes H_t(E') \stackrel{\times}{\to} H_{s+t}(E \times E') \stackrel{\tau_*}{\to} \\ 
\tilde{H}_{s+t}((E\times _M E')^{TM}) \stackrel{u_*}{\to} \tilde{H}_{s+t-d}(E\times _M E') \stackrel{m'_*}{\to} \tilde{H}_{s+t-d}(E').
\end{multline*}
Again, we will see in Section \ref{sec:ringspectra} that this product makes $\mathbb{H}_*(E')$ into a module over $\mathbb{H}_*(E).$

\section{Ring spectra}\label{sec:ringspectra}
Given a smooth embedding of closed manifolds $e: P\to N$, along with a vector bundle or virtual bundle $\zeta \to N$, we obtain a Pontrjagin-Thom collapse map $$\tau: N^{\zeta} \to P^{e^*\zeta\oplus \nu_e}.$$
If $\zeta$ is a virtual bundle then $N^{\zeta}$ is now the Thom spectrum of $\zeta.$  In the case $\zeta=-E$ for a $k$-dimensional vector bundle $E$, $N^{\zeta}$ is defined as $$N^{-E} = \Sigma ^{-(N+k)}N^{E^{\perp}}$$ where $E^{\perp}$ is the $N$-dimensional orthogonal complement bundle to an embedding $E\into N\times \mathbb{R}^{k+N}.$  This generalization of the Pontrjagin-Thom construction also extends to the situation of \cite{ralph}: if 
$$\begin{CD}
E @>f>> X \\
@VqVV @VVpV \\
P @>>e> N
\end{CD}$$
is homotopy Cartesian and $e$ is a smooth embedding of closed manifolds, then for any vector or virtual bundle $\zeta$ over $X$ we obtain a Pontrjagin-Thom collapse map
$$\tau: X^{\zeta} \to E^{f^*\zeta \oplus q^*\nu_e}.$$ 

We will use this construction to prove the following theorem.

\begin{thm}\label{ringsp} Let $F\to E\stackrel{p}{\to} M$ be a smooth fiberwise monoid with $M$ closed of finite dimension $d$.  Then the Thom spectrum $E^{-TM}$ is an associative ring spectrum with unit.  Furthermore, the induced map $p:E^{-TM} \to M^{-TM}$ is a map of ring spectra. \end{thm}

By an associative ring spectrum with unit we mean in a weak sense; that is, it is a spectrum $X$ along with a map $\mu : X \wedge X \to X$ and a map $\iota:S \to X$, where $S$ denotes the sphere spectrum, so that the diagram
$$\begin{CD}
X\wedge X \wedge X @>{\mu \wedge id}>> X\wedge X \\
@V{id \wedge \mu}VV @VV{\mu}V \\
X\wedge X @>{\mu}>> X 
\end{CD}$$
commutes up to homotopy and the compositions 
$$X \stackrel{\cong}{\to} S \wedge X \stackrel{\iota \wedge id}{\to} X\wedge X \stackrel{\mu}{\to} X$$
and 
$$X \stackrel{\cong}{\to} X \wedge S \stackrel{id \wedge \iota}{\to} X\wedge X \stackrel{\mu}{\to} X$$
are homotopic to the identity on $X$.
Note that in this paper, most maps of ring spectra will be well-defined only up to homotopy.  The first-named author intends to rigidify this structure in future work \cite{thesis}.

As before, let $G\into E\stackrel{p}{\to} M$ be a fiberwise monoid with multiplication map $m$.  Assume $M$ is a smooth, closed manifold of dimension $d$.  We saw previously that we can apply the Pontrjagin-Thom construction to  $\tilde{\Delta}: E\times _M E \into E\times E$ with  $\nu_{\tilde{\Delta}}\cong p^*(TM).$  Pulling back the virtual bundle $-TM\times -TM$ over $E\times E$, we obtain a Pontrjagin-Thom map
$$\tau: (E\times E)^{(p\times p)^*(-TM\times -TM)} \to (E\times_M E)^{\tilde{\Delta}^*(p\times p)^*(-TM\times -TM)\oplus p^*(TM)}.$$
In simplified notation this is $$\tau: (E\times E)^{(-TM \times -TM)} \to (E\times _M E)^{-2TM\oplus TM}$$
or, in particular, 
$$\tau:E^{-TM}\wedge E^{-TM} \to (E\times _M E)^{-TM}.$$
Furthermore, by the commutativity of 
$$\begin{CD}
E\times _M E @>m>> E \\
@V{p}VV @VV{p}V \\
M @= M
\end{CD}$$
we see that $m$ induces a map of Thom spectra
$$m: (E\times _M E)^{-TM} \to E^{-TM}.$$
Thus, we define the ring spectrum structure on $E^{-TM}$ as the composition
$$\mu: E^{-TM}\wedge E^{-TM} \stackrel{\tau}{\to} (E\times _M E)^{-TM} \stackrel{m}{\to} E^{-TM}.$$

Notice that we have used two different methods to obtain maps between Thom spectra: Pontrjagin-Thom collapse maps and maps induced from maps of bundles.  In confirming the properties of the ring spectrum $E^{-TM}$ we will often use the following naturality property of the Pontrjagin-Thom maps with respect to the induced maps.
\begin{lemma} \label{nat}
Let $e: P \into N$ be a smooth embedding of closed  manifolds. Suppose we have two fiber bundles $X_1\stackrel{p_1}{\to} N$ and $X_2 \stackrel{p_2}{\to} N$ and a map $g:X_1 \to X_2$ such that $p_2 g = p_1$.  Then we have the pullback diagram
$$\begin{CD}
E_1 := e^*X_1 @>{f_1}>> X_1 \\
@V{g}VV @VV{g}V \\
E_2 := e^*X_2 @>{f_2}>> X_2 \\
@V{q_2}VV @VV{p_1}V \\
P @>{e}>> N 
\end{CD}$$
Let $q_1 = q_2g.$  Let $\zeta$ be a vector or virtual bundle over $X_2.$  Then the following diagram of Thom spectra commutes:
$$\begin{CD}
X_1^{g^*\zeta} @>{\tau}>> E_1^{f_1^*g*\zeta \oplus q_1^*\nu_e} = E_1 ^{g^*(f_2^*\zeta \oplus q_2^*\nu_e)} \\
@VgVV @VVgV \\
X_2^\zeta @>{\tau}>> E_2^{f_2^*\zeta \oplus q_2^*\nu_e}.
\end{CD}$$
\end{lemma}
\begin{proof}
This follows easily from \cite{ralph}, since after replacing $p_1$ and $p_2$ by path fibrations we can use the same tubular neighborhood $\eta_e$ to obtain the collapse maps for both $f_1$ and $f_2$.  
\end{proof} 

To check that the ring spectrum structure on $E^{-TM}$ is associative, we will apply the lemma to
$$\begin{CD}
E\times _M E \times _M E @>{id \times \tilde{\Delta}}>> E\times _M E \times E \\
@V{m\times id}VV @VV{m\times id}V \\
E\times_ M E @>{\tilde{\Delta}}>> E\times E.
\end{CD}$$
Here the embedding $e$ will be $\Delta:M \to M\times M$ and $\zeta = -TM \times- TM.$  Hence the following diagram commutes:
$$\begin{CD}
(E\times _M E)^{-TM}\wedge E^{-TM}  @>{\tau}>> (E\times _M E\times _M E) ^{-TM} \\
@V{m\times id}VV @VV{m\times id}V \\
 E^{-TM}\wedge E^{-TM}  @>{\tau}>> (E\times _M E)^{-TM}.
\end{CD}$$
Hence $\mu\circ (\mu \wedge id)$ is given by $$E^{-TM} \wedge E^{-TM} \wedge E^{-TM} \stackrel{\tau}{\to} (E\times _M E \times _M E)^{-TM} \stackrel{m\circ(m\times id)}{\longrightarrow} E^{-TM}.$$
Since $m\circ (m\times id) = m\circ(id\times m)$, it is now clear that $\mu \circ (\mu \wedge id) = \mu \circ (id \wedge \mu).$  

The Thom spectrum $M^{-TM}$ also has the structure of an associative ring spectrum, given by applying the Pontrjagin-Thom construction to the diagonal embedding $\Delta:M\into M\times M$ and using the virtual bundle $-TM \times -TM,$
$$\tau: M^{-TM} \wedge M^{-TM} \to M^{-TM}.$$  
This ring spectrum structure was shown to rigidify to a commutative symmetric ring spectrum in \cite{cohen}.  After applying the homology Thom isomorphism on both sides, $\tau$ realizes the intersection product on $H_*(M).$  Furthermore, $M^{-TM}$ is S-dual to $M_+$.  The unit $j:S \to M^{-TM}$ is dual to the projection map $M_+ \to S^0.$  We can now define the unit of $E^{-TM}$ as the composition $$\iota: S \stackrel{j}{\to} M^{-TM} \stackrel{s}{\to} E^{-TM}$$ where $s$ is the unit section.   Hence $E^{-TM}$ is an associative ring spectrum with unit; notice that if the fiberwise multiplication $m$ is commutative then so is $E^{-TM}.$

To show that $p:E^{-TM} \to M^{-TM}$ is a map of ring spectra we just need to see that the diagram below commutes.
$$\begin{CD}
(E\times E)^{-TM\times -TM} @>{\tau}>> (E\times _M E)^{-TM} @>m>> E^{-TM} \\
@VV{p\times p}V @VV{p}V @VVpV \\
(M\times M)^{-TM \times -TM} @>{\tau}>> M^{-TM} @= M^{-TM} 
\end{CD}$$
The square on the left commutes from the construction of the collapse map, and the square on the right commutes since $p\circ m = p.$ This proves Theorem~\ref{ringsp}.

In the case when $M$ is oriented, after applying the Thom isomorphism, the ring spectrum multiplication realizes the homology product described in Section~\ref{sec:fibmon}.  That is, 
the following diagram commutes, where $u_*$ denotes the Thom isomorphism:
$$\begin{CD}
H_{q-2d}(E^{-TM}\wedge E^{-TM}) @>{\mu _*}>> H_{q-2d}(E^{-TM}) \\
@V{u_*}V{\cong}V @V{\cong}V{u_*}V \\
H_q(E\times E) @>{\mu_*}>> H_{q-d}(E)
\end{CD}$$
as shown in \cite{cj}.

Hence we have proved the following:

\begin{prop} \label{eme}
Let $p:E \to M$ be a fiberwise monoid over a closed, smooth, oriented $d$-dimensional manifold $M$.  Then $\mathbb{H}_*(E)$ has a graded algebra structure.  Moreover, $p_*:\mathbb{H}_*(E) \to \mathbb{H}_*(M)$ is an algebra homomorphism, where $\mathbb{H}_*(M)$ is equipped with the intersection product. 
\end{prop}

The analogous construction works for fiberwise modules.  If $E' \to M$ is a fiberwise module over a fiberwise monoid $E \to M$, then the Thom spectrum ${E'}^{-TM}$ is a module over the ring spectrum $E^{-TM}$ via
$$E^{-TM} \wedge {E'}^{-TM} \stackrel{\tau}{\to} (E \times _M {E'})^{-TM} \stackrel{m'}{\to} {E'}^{-TM}.$$  If $M$ is oriented, applying the Thom isomorphism in homology at all stages realizes the module structure $\mu_*$ from Section~\ref{sec:fibmon}.

\begin{prop}
Let $E' \to M$ be a fiberwise module over a fiberwise monoid $E\to M.$  Then $\mathbb{H}_*(E')$ is a module over the associative algebra $\mathbb{H}_*(E).$
\end{prop}

Cohen, Jones and Yan proved the following theorem for the case $E=LM$ in \cite{cjy} and their proof easily generalizes for fiberwise monoids.

\begin{thm}\label{specseq}
Let $F\into E \to M$ be a fiberwise monoid with $M$ a smooth, closed, oriented, simply connected manifold.  Then there is a second quadrant spectral sequence of algebras $\{E^r_{p,q}, d^r : p \leq 0, q\geq 0\}$ such that 
\begin{itemize}
\item{$E^r_{*,*}$ is an algebra and the differential $d^r: E^r_{*,*} \to E^r_{*-r,*+r-1}$ is a derivation for each $r\geq 1.$}
\item{The spectral sequence converges to $\mathbb{H}_*(E)$ as algebras.}
\item{For $m,n \geq 0$, 
$$E^2_{-m,n} \cong H^m(M; H_n(F)).$$
Furthermore, $E^2_{-*,*}\cong H^*(M; H_*(F))$ as algebras, where the algebra structure on $H^*(M;H_*(F))$ is given by the cup product on the cohomology of $M$ with cohomology in the Pontrjagin ring $H_*(F).$}
\end{itemize}
\end{thm}

\section{Naturality}\label{sec:nat}
In this section we will prove the following naturality theorem for the ring spectra arising from fiberwise monoids.
\begin{thm}\label{naturality}
\begin{enumerate}
\item\label{morph1} A  morphism of  fiberwise monoids of the form
$$\begin{CD} 
F @= F \\
@VVV @VVV \\
 E_1 = f^*(E_2) @>>> E_2 \\
@VVV @VVV \\ 
M_1 @>f>> M_2 \end{CD} $$
 with $M_i$ closed and $f$ smooth, induces a map of ring spectra  $$\theta_f :E_2^{-TM_2} \to E_1^{-TM_1}$$ that is compatible with the Atiyah dual of f, $$M_2^{-TM_2} \to M_1^{-TM_1}.$$  

\item\label{morph2} A morphism of  fiberwise monoids of the form
$$\begin{CD}
F_1 @>f>> F_2 \\
@VVV @VVV \\
E_1 @>f>> E_2 \\
@VVV @VVV \\
M @= M
\end{CD}$$
with $M$ closed, induces a map of ring spectra $$f: E_1^{-TM} \to E_2^{-TM}$$ that is compatible with the identity map of $M^{-TM}.$
\end{enumerate}
\end{thm}

\begin{proof}
Suppose we have a morphism of fiberwise monoids as in (\ref{morph1}) above.  
  Fix a smooth embedding $e$ of $M_1$ into a sphere $S^N$ of sufficiently large dimension.  Then 
$$\begin{CD} E_1 @>{f\times (e\circ p_1)}>> E_2 \times S^N \\
@V{p_1}VV @VV{p_2 \times id}V \\
M_1 @>{f\times e}>> M_2 \times S^N
\end{CD}$$
is also a morphism of fiberwise monoids, and $E_1 = (f\times e)^*(E_2 \times S^N).$  Hence there is a Pontryagin-Thom collapse map 
$$E_2^{-TM_2}\wedge {S^N}^{-TS^N} \to E_1^{(f\times e)^*(-TM_2 \times -TS^N) \oplus \nu_{(f\times e)}} = E_1^{-TM_1}$$ which we will denote by $\tau _f.$  Now ${S^N}^{-TS^N}$ is ring spectrum with unit $j:S \to {S^N}^{-TS^N}.$  Hence we have a map of spectra
$$\theta _f :E_2^{-TM_2} \cong E_2^{-TM_2}\wedge S \stackrel{id \wedge j}{\to} E_2^{-TM_2}\wedge {S^N}^{-TS^N} \stackrel{\tau _f}{\to} E_1^{-TM_1}.$$   

We now verify that $\theta _f$ is a map of ring spectra.  First consider the diagram
$$\begin{CD}
E_2^{-TM_2}\wedge E_2 ^{-TM_2} \wedge S \wedge S @>{id\wedge id \wedge j \wedge j}>> E_2^{-TM_2} \wedge E_2^{-TM_2} \wedge {S^N}^{-TS^N} \wedge {S^N}^{-TS^N} \\
@V{\mu _2 \wedge \alpha}VV @VV{\mu _2 \wedge \tau_\Delta}V \\
E_2^{-TM_2} \wedge S @>{id \wedge j}>> E_2^{-TM_2}\wedge {S^N}^{-TS^N}.
\end{CD}$$
Here $\alpha$ is the identification $S\wedge S = S$ (which is the multiplication on the ring spectrum $S$) and $\tau_\Delta$ is the Pontrjagin-Thom collapse map from $\Delta: S^N \into S^N \times S^N$, hence is the multiplication for the ring spectrum ${S^N}^{-TS^N}.$  This diagram commutes since $j$ is the unit of ${S^N}^{-TS^N}.$
Next consider
$$\begin{CD} 
E_2^{-TM_2} \wedge  E_2^{-TM_2} \wedge {S^N}^{-TS^N} \wedge {S^N}^{-TS^N} @>{\tau _f \wedge \tau _f}>> E_1^{-TM_1}\wedge E_1^{-TM_1} \\
 @V{\tau_{\tilde{\Delta}} \wedge \tau_\Delta}VV @VV{\tau_{\tilde{\Delta}}}V \\
 (E_2 \times _{M_2} E_2)^{-TM_2} \wedge {S^N}^{-TS^N} @>{\tau_g}>> (E_1 \times _{M_1} E_1)^{-TM_1} \\
 @V{m_2 \wedge id}VV @VV{m_1}V \\
 E_2^{-TM_2} \wedge {S^N}^{-TS^N} @>{\tau _f}>> E_1^{-TM_1}
\end{CD}$$
where $g = f \times (e\circ p_1): E_1\times _{M_1} E_1 \into E_2 \times _{M_2} E_2 \times S^N.$  

 The commutativity of 
$$\begin{CD} E_2 \times E_2 \times S^N \times S^N @<{\hookleftarrow}<< E_1 \times E_1 \\
@A{\tilde{\Delta} \times \Delta}AA @AA{\tilde{\Delta}}A \\
E_2 \times _{M_2} E_2 \times S^N @<{\hookleftarrow}<< E_1 \times _{M_1} E_1
\end{CD}$$
and the naturality of the Pontrjagin-Thom construction show that the top square commutes.  To check the commutativity of the bottom square, we apply the lemma to
$$\begin{CD} E_1\times _{M_1} E_1 @>{\into}>> E_2\times _{M_2} E_2 \times S^N \\
@V{m_1}VV @VV{m_2 \times id}V \\
E_1 @>{\into}>> E_2 \times S^N
\end{CD} $$
with $P=M_1$, $N = M_2 \times S^N$ and $\zeta = -TM_2 \times -TS^N$.  
Combining the two diagrams above shows that 
$$\begin{CD}
E_2^{-TM_2}\wedge E_2^{-TM_2} @>{\theta _f \wedge \theta _f}>> E_1^{-TM_1} \wedge E_1^{-TM_1} \\
@V{\mu _2}VV @VV{\mu _1}V \\
E_2^{-TM_2} @>{\theta _f}>> E_1^{-TM_1}
\end{CD}$$
commutes.

We also check that $\theta _f$ takes the unit of $E_2^{-TM_2}$ to the unit of $E_1^{-TM_1}.$  In the diagram
$$\begin{CD}
S^0 @>{j}>> M_2^{-TM_2} @>{s_2}>> E_2^{-TM_2} \\
@| @VV{\theta _f}V @VV{\theta _f}V \\
S^0 @>{j}>> M_1^{-TM_1} @>{s_1}>> E_1^{-TM_1}
\end{CD}$$
the first square commutes because it is dual to
$$\begin{CD}
S^0 @<<< {M_2} _+ \\
@| @AAA \\
S^0 @<<< {M_1} _+
\end{CD} $$
which commutes. The second square commutes by the construction of the collapse maps.   This proves the first part of Theorem~\ref{naturality}. 

To check the second part, suppose we have a  morphism of fiberwise monoids
$$\begin{CD}
F_1 @>>> F_2 \\
@VVV @VVV\\
E_1 @>f>> E_2 \\
@V{p_1}VV @VV{p_2}V \\
M @= M.
\end{CD}$$
Then $f$ induces a map of Thom spectra
$$f: E_1 ^{-TM} \to E_2 ^{-TM}.$$
This map is a map of ring spectra since in the diagram
$$\begin{CD}
E_1 ^{-TM} \wedge E_1 ^{-TM} @>{\tau}>> (E_1 \times _M E_1)^{-TM} @>{m_1}>> E_1 ^{-TM} \\
@V{f\wedge f}VV @VVfV @VVfV \\
E_2^{-TM} \wedge E_2^{-TM} @>{\tau}>> (E_2 \times _M E_2)^{-TM} @>{m_2}>> E_2^{-TM}
\end{CD}$$
both squares commute - the first by the lemma, and the second since $f\circ m_1 = m_2 \circ f.$  It is clear that $f$ sends the unit of $E_1^{-TM}$ to the unit of $E_2^{-TM}$, since $f\circ s_1 = s_2.$
\end{proof}

The unit section induces a map of the trivial fiberwise monoid over $M$ into $E$ of the form in part (ii) of the theorem, hence we obtain: 
\begin{cor} The unit section $s: M \to E$ of a fiberwise monoid over a closed, smooth, oriented $d$-dimensional manifold induces an algebra homomorphism $s_*: \mathbb{H}_*(M) \to \mathbb{H}_*(E)$.
\end{cor}

 If we apply part (i) of the theorem 
 to the case when $f$ is the inclusion of the base point into $M$ and $E_2=E$ then we get a map of ring spectra
 $E^{-TM} \to \Sigma^{\infty}F_+$.  
This induces in homology a homomorphism $c: \mathbb{H}_*(E) \to H_*(F)$ called in \cite{cs} transverse intersection with a fiber.
Thus we obtain:
\begin{prop} \label{ef}
The intersection with a fiber is a homomorphism of algebras $c: \mathbb{H}_*(E) \to H_*(F)$, where $H_*(F)$ is equipped with the Pontrjagin product.
\end{prop}

Likewise we obtain naturality properties for fiberwise modules.
\begin{prop}\label{naturality2}
\begin{enumerate}
\item
Let $E' \to M_2$ be a fiberwise module over a fiberwise monoid $E \to M_2$, and let $f:M_1 \to M_2$ be a smooth map of closed, oriented manifolds.  Then $f^*E' \to M_1$ is a fiberwise module over the fiberwise monoid $f^*E \to M_1$ and the diagram 
$$\begin{CD}
E^{-TM_2} \wedge {E'}^{-TM_2} @>{\mu '}>> {E'}^{-TM_2} \\
@V{\tau_f \wedge \tau_f}VV @VV{\tau_f}V \\
f^*E^{-TM_1} \wedge f^*{E'}^{-TM_1} @>{\mu '}>> f^*{E'}^{-TM_1}
\end{CD}$$
commutes.
\item
Let $E' \to M$ and $E'' \to M$ be fiberwise modules over the smooth fiberwise monoid $E\to M$, with a smooth module map
$$\begin{CD}
E' @>g>> E'' \\
@VVV @VVV \\
M @= M.
\end{CD}$$
Then $g:{E'}^{-TM} \to {E''}^{-TM}$ is a map of $E^{-TM}$ modules.
\end{enumerate}
\end{prop}

\section{Little cubes and string topology} \label{sec:c_n}

In this section we study the operations induced by the little
$n$-cubes operad in a fiberwise setting.

We recall that a topological operad $C$ is a sequence of spaces
$C(n)$, together with an action of $\Sigma_n$ on $C(n)$, a unit in $C(1)$, 
and structure maps
$$C(k) \x C(i_1) \x \dots \x C(i_k) \to C(i_1 +\dots+i_k)$$
satisfying appropriate associativity, equivariance and unit axioms
\cite{gk}.

The little $n$-cubes operad $\CC_n$ is defined so that
$\CC_n(k)$ is the space of embeddings of the disjointed union
of $k$ copies of a cube
into a cube  $\coprod_k I^n \to I^n$,
that on each copy is a product of affine embeddings.
The operad structure maps are defined by composition and the symmetric group action
 by reordering.
By convention $\CC_n(0)$ is a point.

A space $A$ is an {\em algebra} over a topological operad $C$ 
if there are structure maps
$\theta_n:C(k) \x A^k \to A$ satisfying appropriate associativity, equivariance and unit
conditions \cite{gk}.

A $\CC_n$-algebra is usually also called a $C_n$-{\em space}.
A typical example of a $\CC_n$-algebra is an $n$-fold loop space
$\Om^n(X)$, the space of based maps from $S^n$ to a topological space $X$
(Section \ref{examples}).

\

Operads and their algebras can be defined similarly in any symmetric
 monoidal category, if we replace
the cartesian product above by the tensor product.

We are interested in the category of spaces over a fixed space $M$.
Such category is symmetric monoidal
by its categorical product: the tensor product of two objects
$X \to M$ and $Y \to M$ is the pullback $X\x_M Y$ equipped with the obvious
map to $M$.
Given a topological operad $O$,
the collection of trivial bundles $O(k) \times M \to M$ is an operad $\bar{O}$
in the category of spaces over $M$.
\bde
We say that a bundle $\pi:E \to M$ is a  fiberwise $O$-algebra
if it is an algebra over $\bar{O}$ in the category of spaces over $M$.
\ede

In particular each fiber $\pi^{-1}(m)$
over a point $m \in M$ is an algebra over $O$.
If $O(0)$ is a point,
this yields a section $s:M \to E$ of $\pi$.
If $O$ is the operad defining monoids, then a fiberwise $O$-algebra
is a fiberwise
monoid as introduced earlier.

If $O=\CC_n$ we say also that $\pi:E \to M$
is a fiberwise $C_n$-space.
A typical example of fiberwise $C_n$-space is
the free mapping space $map(S^n,X)$ equipped with the evaluation
projection to $X$ (Section \ref{examples}).

Another monoidal category of interest for us is the category of spectra equipped
with the smash product.
The suspension spectrum functor turns $\CC_n$ into an operad
$\Sigma^{\infty}(\CC_n)_+$ in the category of spectra. An algebra
over such operad is called an {\em $E_n$-ring} spectrum.
In the rest of the paper we will mean the weak version of this notion
in the sense that the
 associativity axiom is relaxed up to homotopy.

We are ready now to state the main result of the section.

\begin{thm}      \label{e_n}

Let $M$ be a smooth, closed, compact manifold.
Let $\pi: E \to M$ be a fiber bundle that is a fiberwise $C_n$-space.
Then $E^{-TM}$ is an $E_n$-ring spectrum.

\end{thm}

\begin{proof}
As done earlier we write $-TM$ also for its pullback via $\pi$ or other
projections to $M$.
By definition we have structure maps
$$\theta_k:C_n(k) \x (\pi^k)^{-1}(\Delta_k(M)) \to E,$$
where $M \cong \Delta_k(M) \subset M^k $ is the thin diagonal.
On the other hand we have a pullback diagram

$$\begin{CD}
 (\pi^k)^{-1}(\Delta_k(M)) @>>>  E^k  \\
@VVV                                    @VV\pi^kV             \\
 \Delta_k(M)             @>i_k>>         M^k \,.
\end{CD}$$

The embedding $i_k$ has as normal bundle $\nu_k$, the kernel of the projection
$\oplus_k TM \to TM$, that is the reduced representation of $\Sigma_k$ on
$TM$, and is
isomorphic to $\oplus_{k-1}TM$ as a non-equivariant vector bundle.
On the level of $\Si_k$-equivariant virtual bundles the equality
$\oplus_k (-TM) \oplus \nu_k \cong -TM$ holds.
Thus the Thom-Pontrjagin construction in \cite{ralph}
twisted by the virtual bundle $(-TM)^k$ on $E^k$ yields  
 a map of spectra
$$\tau_k:(E^{-TM})^{\sm k} \to (\pi^k)^{-1}(\Delta_k(M))^{-TM}.$$
The structure map of $E$ gives us on spectra
$$\theta_k^{-TM}:\Sigma^{\infty}C_n(k)_+ \sm (\pi^k)^{-1}(\Delta_k(M))^{-TM}
\to E^{-TM}.$$
 The desired structure map of $E^{-TM}$ is then the composite
$$\theta_k^{-TM} \circ ( \Sigma^{\infty}C_n(k)_+  \sm \tau_k):
\Sigma^{\infty}C_n(k)_+ \sm (E^{-TM})^{\sm k} \to E^{-TM}.$$

\

The equivariance is immediate. The associativity follows 
from Lemma \ref{nat} similarly as for monoids, using the associativity
condition for $E$. The unit is constructed just as in the case of monoids.

\end{proof}

\begin{prop} \label{en-maps}
The projection $p:E \to M$ and the section $s:M \to E$ induce
maps of $E_n$-ring spectra $E^{-TM} \to M^{-TM}$ and $M^{-TM} \to E^{-TM}$.
The inclusion of the fiber $F \to E$ induces a map of $E_n$-ring
spectra $E^{-TM} \to \Sigma^{\infty}(F_+)$.
\end{prop}

\begin{proof}
 Apply Theorem \ref{naturality} (1) to the projection and the section,
and Theorem \ref{naturality} (2) to the fiber inclusion.
\end{proof}

\

It is well known that $C_n$-spaces and $E_n$-ring spectra have a
wealth of homology operations, thoroughly studied by F. Cohen in
\cite{lnm533}.

As a consequence of Theorem \ref{e_n} we obtain :

\begin{cor}
Let $E \to M$ be a fiberwise $C_n$-space
 over a closed, compact,
smooth manifold $M$. We assume $M$ to be oriented unless we work at the prime 2.

There exist operations

\

{\rm (Prime 2)} $Q_i:\hh_q(E,\Z_2) \to \hh_{2q+i}(E,\Z_2)$, for $0 \leq i \leq n-1$;

\

{\rm (Odd primes)} $\beta^\varepsilon Q_i:\hh_q(E,\Z_p) \to \hh_{pq+i(p-1)-\varepsilon}(E,\Z_p)$,
 for $p$ odd prime, $q+i$ even, $0 \leq i \leq n-1, \,\ep =0,1;$ 

\

{\rm (Browder)}  
$\lambda_{n-1}: \hh_p(E) \ot \hh_q(E) \to \hh_{p+q+n-1}(E)$.

\

The $Q_i$'s operations satisfy the Adem relations and the internal Cartan relations (including unstable relation Thm 1.3(2) in \cite{lnm533}).
The operation $\lambda_{n-1}$ satisfies the internal Cartan formula, the
 Jacobi identity,
and the compatibility with $Q_i$'s in Thm 1.2 (8), Thm 1.3 (4,5) of \cite{lnm533}.

\end{cor}
\begin{proof}
These operations are always defined on the homology of
 an $E_n$-ring spectrum.
Namely $\lambda_{n-1}(x \ot y) = (\theta_2)_*(e \ot x \ot y)$,
with $e \in H_{n-1}(\CC_n(2)) \cong H_{n-1}(S^{n-1})$ a generator
and $\theta_2:\Sigma^{\infty} \CC_{n}(2)_+ \sm E^{-TM} \sm E^{-TM} \to E^{-TM}$ the structure map.
  The operations $Q_i's$ depend on the choice of classes
$e_{i(p-1)-\varepsilon} \in H_{i(p-1)-\varepsilon}(\CC_n(p)/\Sigma_p,\pm \Z_p)$  (untwisted or twisted by
the sign representation) 
 and are defined by
$Q_i(x) = (\bar{\theta}_p)_*(e_{i(p-1)-\varepsilon} \ot x^{\ot p})$, where
$\bar{\theta}_p:\Sigma^{\infty}\CC_n(p)_+ \sm_{\Sigma_p} (E^{-TM})^{\sm p} \to
E^{-TM}$ is induced by the structure map,
  and we use the fact that $$H_*(\Sigma^{\infty}\CC_n(p)_+ \sm_{\Sigma_p} (E^{-TM})^{\sm p})
\cong \mathcal{H}_*(\CC_n(p)/\Sigma_p, H_*(E^{-TM})^{\ot p}).$$
 Via the Thom isomorphism 
$\hh_*(E) \cong H_*(E^{-TM})$ this defines our operations.
 All relations hold at the level of equivariant homology of the
operad of little $n$-cubes.

\end{proof}

\

{\em Remark}:
In the special  case $E=map(S^n,M)$
 the homological action of the $n$-little cubes operad extends
to an action of the $(n+1)$-little cubes operad. See \cite{cj} for $n=1$ and \cite{Hu} for
$n>1$. In this case there are additional homology operations $Q_n$, and $\lambda_{n-1}$ is trivial. 
They were constructed in \cite{Tourtchine} and \cite{Westerland}.

\

In the special case $E=M$ we have that $M^{-TM}$ is an $E_\infty$-ring
spectrum, and this defines the operations above on $\hh_*(M)$.
We will show that they are Poincar\'e dual to Steenrod operations.

Let us denote the homology operation Poincar\'e dual to  
the Steenrod square $Sq^i:H^q(M,\Z_2) \to H^{q+i}(M,\Z_2)$
by
$PD(Sq^i):\hh_{-q}(M,\Z_2) \to \hh_{-q-i}(M,\Z_2)$, with our shift convention.
Similarly the Poincar\'e
dual of the power operation $P^i$ at the prime $p$ is denoted
 $$PD(P^i):\hh_{-q}(M,\Z_p) \to \hh_{-q-2i(p-1)}(M,\Z_p).$$

\begin{thm}
Let $M$ be a closed smooth manifold.

(i) For $x \in \hh_{-q}(M,\Z_2)$,
$$Q_i(x) =PD(Sq^{q-i})(x) \in \hh_{-2q+i}(M,\Z_2).$$

(ii) For $M$ oriented and $x \in \hh_{-q}(M,\Z_p)$, with $p$ odd,
$$\beta^\varepsilon Q_i(x)=PD(\beta^\varepsilon P^{\frac{q-i}{2}})(x)
\in \hh_{pq+i(p-1)}(M,\Z_p).$$

(iii) For $M$ oriented the operations $\lambda_n$ vanish for any $n$.

\end{thm}
\begin{proof}
We give a proof of (i).
Recall that the Serre spectral sequence of
$$M^2 \to S^i \x_{\Z_2} M^2 \to \RP^i$$
 collapses at the prime 2.
By definition $Q_i(x)=\pi_* j_!([\RP^i] \ot x \ot x)$, where
$j:\RP^i \x M \to S^i \x_{\Z_2} M^2$ is the embedding
 induced by the diagonal,
 $j_!$ is induced by the Thom-Pontrjagin
 collapse followed by Thom isomorphism,
 and
  $\pi:\RP^i \x M \to M$ is the projection.
It is well known that $j_!$ in homology corresponds under Poincar\'e duality
to $j^*$
in cohomology.

Let $x' \in H^q(X,\Z_2)$ be the cohomology class dual to $x$.
Then $$[\RP^0]^*\ot x' \ot x' \in H^{2q}(S^i \x_{\Z_2} M^2,\Z_2)$$
is dual to $[\RP^i] \ot x \ot x$.
Thus $j_!([\RP^i] \ot x \ot x)$ is Poincar\'e dual to
$j^*([\RP^0]^*\ot x' \ot x')= \sum_{t=0}^i [\RP^t]^* \ot Sq^{q-t}(x')$, namely it is
$\sum_{t=0}^i [\RP^{i-t}]\ot PD(Sq^{q-t})(x)$, but the projection $\pi_*$
kills all summands except the last.

The proof of (ii) is very similar, and is done with respect
to the diagonal map $S^{2N+1}/\Z_p \x M \to S^{2N+1} \x_{\Z_p}  M^p$,
a finite dimensional approximation of
$B\Z_p \x M \to E\Z_p \x_{\Z_p} M^p$.

Statement (iii) follows because the operation $\lambda_n$
 is an obstruction to
the extension from  an $E_{n+1}-$ to an $E_{n+2}$-ring structure on $M^{-TM}$.
\end{proof}

Statements (i),(ii) in
the following proposition are similar to Propositions \ref{eme}, \ref{ef}
and follow from Proposition \ref{en-maps}.

\bpr
Let $\pi:E \to M$ be a fiberwise $C_n$-space over an oriented  manifold 
$M$, with fiber $F$ and canonical section $s$. Then

(i) The homomorphisms $\pi_*:\hh_*(E) \to \hh_*(M)$ and $s_*:\hh_*(M) \to \hh_*(E)$
commute with the homology operations $Q_i's$ and $\lambda_{n-1}$.

(ii) Let $c:\hh_*(E) \to H_*(F)$ be the transverse intersection
with a fiber.
Let us denote by $Q_i$ and $\lambda_{n-1}$ the operations
of $C_n$-space on $H_*(F)$ \cite{lnm533}.
Then the homomorphism $c$ commutes with the operations $Q_i's$ and
$\lambda_{n-1}$.
\epr

\

There is a result similar to Theorem \ref{e_n} for modules.
Let $O$ be an operad in a fixed
symmetric monoidal category $\mathcal{C}$.
 A module $M$ over an $O$-algebra $A$
is an object of $\mathcal{C}$
 equipped with structure maps
$O(k)\ot A^{\ot k-1} \ot M \to M$ satisfying appropriate axioms
\cite{gk}.
We consider modules in the category of spectra in the weak sense,
requiring that the appropriate diagrams in the axioms commute up to homotopy,
as we did for algebras.

\begin{defn}
A bundle
 $E' \to M$ is a fiberwise module over a fiberwise $O$-algebra $E \to M$
if $E' \to M$ is a module over the $\bar{O}$-algebra $E \to M$ in the
category of spaces over $M$.
\end{defn}

A typical example of fiberwise module over a fiberwise $C_n$-space
is the mapping space $map(N,X)$ from a $n$-dimensional manifold,
that is a module over $map(S^n,X)$ (section \ref{examples}).

\begin{prop}
Let $E' \to M$ be a fiberwise module over a fiberwise $C_n$-space $E \to M$,
with $M$ a closed compact manifold.

Then $E'^{-TM}$ is a module over the $E_n$-ring spectrum $E^{-TM}$.
\end{prop}

We recall that a $n$-algebra (see \cite{SW}) is an algebra over the operad $H_*(\CC_n,\Q)$.
A 2-algebra is also called a Gerstenhaber algebra.

\begin{cor}
If a bundle $E \to M$ is a fiberwise $C_n$-space and a bundle
 $E' \to M$ is a module over it, 
then $\hh_*(E',\Q)$ is a module
over the $n$-algebra $\hh_*(E,\Q)$.
\end{cor}

\section{ String Topology of Classifying Spaces}\label{sec:bg}
The naturality properties of fiberwise monoids allow us to define and study the string topology of $BG$, when $G$ is a compact Lie group, which is not covered by Chas-Sullivan or Cohen-Jones since $BG$ is not a closed manifold.  We will first study fiberwise monoids over colimits of manifolds.  
Suppose $F \into E \to X$ is a fiberwise monoid with multiplication $m$ and unit section $s$, where $X$ is a direct limit of  smooth, closed, oriented manifolds $\{ X_{\alpha} \}$, 
$$X = \colim_{\alpha \in D} X_\alpha$$
where $D$ is a directed set.  Let $i_{\alpha}$ be the map $X_{\alpha} \to X$, and $i_{\alpha \beta}$ the map $X_\alpha \to X_\beta$ whenever $\alpha \leq \beta$.  Then
 $$E_\alpha := i_{\alpha}^*(E) = \{(x,v)| x \in X_{\alpha}, v\in E, p(v) = i_{\alpha}(x) \}$$ 
is a fiberwise monoid with multiplication
$$m_\alpha ( (x,v),(x,w) ) = (x, m(v,w))$$ 
and unit section
$$s(x) = (x, s(i_\alpha (x))).$$  
Then $E_{\alpha} ^{-TX_\alpha}$ is a ring spectrum with unit, and by Theorem~\ref{naturality} whenever $\alpha \leq \beta$ we have a map of ring spectra
$$\theta_{i_{ \alpha \beta }}: E_\beta ^{-TX_\beta} \to E_\alpha ^{-TX_\alpha}.$$
Hence we can associate to $E\to X$ the pro-ring spectrum $\{ E_\alpha ^{-TM_\alpha}, \theta_{i_{\alpha \beta}} \}.$ 

In particular, let $G$ be a compact Lie group.  Then there is always a model of $BG$ which is a colimit of closed submanifolds (for instance, there is always a faithful finite dimensional representation of $G$ on a vector space $V$;
 then take as a model of $EG$ the space
  of linear embeddings of $V$ in $\mathbb{R}^\infty$, and take the finite spaces to be linear embedding of $V$ into $\mathbb{R}^n$).
    Furthermore, $LBG$ is homotopy equivalent (in fact, fiberwise homotopy equivalent where the homotopy equivalences are morphism of fiberwise monoids over $BG$) to $Ad(EG)$, where $EG$ is a universal $G$-bundle over $BG$.
\begin{defn}Let $G$ be a compact Lie group and let $BG$ be a model for the classifying space of $G$ which is a colimit of smooth manifolds,
 $$BG = \colim_{n \to \infty} M_n.$$  Assume the connectivity of the inclusion $M_n \into BG$ is a non-decreasing function of $n$ which tends to infinity.
Define the string topology of $BG$ to be the pro-ring spectrum 
$$LBG^{-TBG} := \{E_n^{-TM_n}, \theta_{i_{nm}}\}$$
where $E_n = Ad(i_n^*(EG)).$ 
\end{defn}
 Notice that $Ad(i_n^*(EG)) = i_n^*(Ad(EG)).$ 

This construction yields a pro-ring structure on the pro-homology
$$\dots \mathbb{H}_*(E_{n+1}) \to \mathbb{H}_*(E_n) \to \dots \to \mathbb{H}_*(E_1)$$ where we take coefficients in $\mathbb{Z}/2\mathbb{Z}$ if the $M_n$'s are not appropriately oriented.

\begin{lemma}
The homotopy type of the string topology of $BG$ is well-defined.
\end{lemma}
\begin{proof}
Suppose we have two different filtrations of $BG$ by closed submanifolds,
$$BG = \colim _{n \to \infty} M_n$$ 
and
$$BG = \colim _{m \to \infty} N_m.$$
Define $i_n : M_n \into BG$, $j_m: N_m \into BG$, $E_n = Ad(i_n^*EG)$ and $F_m = Ad(j_m^*EG).$  Assume that $i_n$ is $\phi(n)$ connected and $j_m$ is $\psi(m)$ connected, where $\phi$ and $\psi$ are non-decreasing functions that tend to infinity.  We need to show that the pro-ring spectra
$$ \dots \to E_{n+1}^{-TM_{n+1}} \to E_n^{-TM_n} \to \dots  \to E_1 ^{-TM_1}$$
and 
$$ \dots \to F_{n+1}^{-TN_{n+1}} \to F_n^{-TN_n} \to \dots  \to F_1 ^{-TN_1}$$
have the same homotopy type.  For $n \in \mathbb{Z}$, let $m(n)$ be the smallest $m$ such that $\psi(m) \geq \phi(n).$  Then there exists $f_n: M_n \to N_{m(n)}$ which is $\phi(n)$ connected and so that $j_{m(n)} \circ f_n \simeq i_n.$  There is an induced ring spectrum map $\theta _{f_n}: F_{m(n)}^{-TN_{m(n)}} \to E_n ^{-TM_n}.$  In homology, after applying the Thom Isomorphism on both sides, $\theta _{f_n}$ is $({f_n})_!$.  So $\theta _{f_n}$ is an isomorphism on $H_*$ for $g\geq * \geq g - \phi(n)$, where $g=dim(G).$  It is now clear that the maps $\{ f_n \}$ give a homotopy equivalence of the two pro-ring spectra.   
\end{proof}

Let us now consider two examples, one with trivial string topology and one with nontrivial string topology.

\medskip
{\bf Example: $G=S^1.$ }
Consider the universal principal $S^1$ bundle $S^1 \into S^\infty \to \mathbb{CP}^\infty.$  Since $S^1$ is abelian, its adjoint bundle is trivial, $Ad(ES^1) = S^1 \times \mathbb{CP}^\infty.$
Now $\mathbb{CP}^\infty = \colim_{n\to \infty} \mathbb{CP}^n.$  Hence the pro-ring spectrum associated to $S^1\times \mathbb{CP}^\infty \to \mathbb{CP}^\infty$ is 
$$\cdots \to (\mathbb{CP}^{n+1} \times S^1)^{-T\mathbb{CP}^{n+1}} \to  (\mathbb{CP}^{n} \times S^1)^{-T\mathbb{CP}^{n}} \to \cdots \to (\mathbb{CP}^{1} \times S^1)^{-T\mathbb{CP}^{1}}.$$
Notice $(\mathbb{CP}^n \times S^1)^{-T\mathbb{CP}^n}= (\mathbb{CP}^n)^{-T\mathbb{CP}^n} \wedge S^1_+.$
By the Thom isomorphism, $H_*((\mathbb{CP}^n \times S^1)^{-T\mathbb{CP}^n}) \cong H_{*+2n}(\mathbb{CP}^n \times S^1)$, so $H_*((\mathbb{CP}^n \times S^1)^{-T\mathbb{CP}^n}) \cong \mathbb{Z}$ for $-2n\leq * \leq 1$ and is $0$ otherwise.  Now the product $$\mu_*: H_s(S^1\times \mathbb{CP}^n) \otimes H_t(S^1 \times \mathbb{CP}^n) \to H_{s+t-2n}(S^1\times \mathbb{CP}^n)$$ is the tensor product of the intersection product on $H_*(\mathbb{CP}^n)$ with the Pontrjagin product on $H_*(S^1)$, on $H_*(S^1 \times \mathbb{CP}^n) \cong H_*(S^1)\otimes H_*(\mathbb{CP}^n).$  As previously noted, this product, with an appropriate dimension shift, describes the product structure on $H_*((\mathbb{CP}^n \times S^1)^{-T\mathbb{CP}^n}).$  In particular, 
$$H_*((\mathbb{CP}^n \times S^1)^{-T\mathbb{CP}^n}) \cong \Lambda(t) \otimes \mathbb{Z}[c]/c^{(n+1)}$$ as rings, where $c$ has degree $-2$ and $t$ has degree $1.$  Furthermore the map $H_*((\mathbb{CP}^{n+1} \times S^1)^{-T\mathbb{CP}^{n+1}}) \to H_*((\mathbb{CP}^n \times S^1)^{-T\mathbb{CP}^n})$ is the usual projection
$$\Lambda(t) \otimes \mathbb{Z}[c]/c^{(n+2)} \to \Lambda(t) \otimes \mathbb{Z}[c]/c^{(n+1)}.$$  The inverse limit of this inverse system of rings is $\Lambda(t) \otimes \mathbb{Z}[[c]].$

\medskip
{\bf Example: $G = O(2).$}
We will use the model $BO(2) = Gr_{2, \infty}$, the Grassmannian of $2$-planes in $\mathbb{R}^\infty$.  Then $BO(2)$ has the obvious filtration by manifolds $M_n = Gr_{2,n}$, the Grassmannian of $2$-planes in $\mathbb{R}^n$, for $n\geq 3.$  As above, let $E_n = i_n^*(EO(2)).$
We will show that $BO(2)$ has nontrivial string topology, in the sense that the pro-ring spectra
$$\dots \to Ad(E_{n+1})^{-TM_{n+1}} \to Ad(E_n)^{-TM_n} \to \dots \to Ad(E_3)^{-TM_3}$$
and 
$$\dots \to (M_{n+1} \times O(2))^{-TM{n+1}} \to (M_n \times O(2))^{-TM_n} \to \dots \to (M_3 \times O(2))^{-TM_3}$$
are not homotopy equivalent.  

Since $M_n$ is not always orientable, $\mathbb{Z}_2$ coefficients will be assumed throughout this example.  Let $\mathcal{A}$ denote the mod $2$ Steenrod algebra.  The $\mathcal{A}$ module structure on $H^*(Ad(E_n)^{-TM_n})$ induces a ``twisted'' $\mathcal{A}$ module structure on $H^*(Ad(E_n))$ defined by $Sq^i_t (x) = (u^*)^{-1}(Sq^i(u^*x))$, where $u^*$ is the Thom isomorphism.  Hence 
$$Sq^i_t(x) = \sum\nolimits_{j=0}^{i}Sq^j(x) \cup w_{i-j}(-TM_n).$$ 
 We then obtain a right $\mathcal{A}$ module structure on $H_*(Ad(E_n))$  defined by 
$$<x,ySq^i_t> = <Sq^i_t x,y>$$
for $x\in H^*(Ad(E_n))$ and $y\in H_{*+i}(Ad(E_n)).$  Notice that $Sq^i_t$ lowers degree by $i$ in homology.  We obviously have this same right $\mathcal{A}$ module structure on $\mathbb{H}_*(Ad(E_n))$.  Furthermore, maps of ring spectra preserve the action of $\mathcal{A}$ on (co)homology.  Thus the inverse system of rings 
$$\dots \to \mathbb{H}_*(Ad(E_{n+1})) \to \mathbb{H}_*(Ad(E_n)) \to \dots \to \mathbb{H}_*(Ad(E_3))$$ 
is also an inverse system of right $\mathcal{A}$ modules.
The $\mathcal{A}$ module structure interacts with the ring multiplication $\mu_*$ as follows:
$$(\mu_*(x,y))Sq^i_t = \mu_*(x\times y)Sq^i_t = \mu_*(\sum_{j=0}^i xSq^j_t \times ySq^{i-j}_t) = \sum_{j=0}^i \mu_*(xSq^j_t,ySq^{i-j}_t)$$
where $x\times y \in H_*(Ad(E_n) \times Ad(E_n))$ and the twisted $\mathcal{A}$ module structure in the second term is induced from $(Ad(E_n) \times Ad(E_n)) ^{-TM \times -TM}$.

Clearly we have the analogous constructions for $\mathbb{H}_*(M_n \times O(2))$, induced from $(M_n \times O(2))^{-TM_n}.$  We will show that the actions of $Sq^1_t$ on $\mathbb{H}_1(Ad(E_n))$ and $\mathbb{H}_1(M_n\times O(2))$ are different.  Notice that $Gr_{2,n}$ is orientable if and only if $n$ is even.  Furthermore, $H^*(Gr_{2,n}) = \mathbb{Z}_2 [w_1, w_2]/R$ where if $\gamma$ is the canonical $2$-plane bundle over $Gr_{2,n}$, $w_i = w_i (\gamma)$ and $R$ is the ideal generated by the relations forced from 
$$(1+w_1 + w_2)(1+w_1(\gamma^\perp) + \dots + w_{n-2}(\gamma^\perp)) = 1.$$  

If $n$ is even, $w_1(-TM_n)=0$ and so the action by $Sq^1_t$ on $\mathbb{H}_1(M_n\times O(2))$ is the same as that by $Sq^1$.  Since $M_n \times O(2)$ is orientable, $xSq^1_t = xSq^1 = 0$ for any $x \in \mathbb{H}_1(M_n\times O(2)).$  

If $n$ is odd, $w_1(-TM_n) = w_1$ and $M_n \times O(2)$ is non-orientable.  $\mathbb{H}^0(M_n \times O(2)) \cong \mathbb{Z}_2 \times  \mathbb{Z}_2 \times  \mathbb{Z}_2 \times  \mathbb{Z}_2 $ with generators $z\times a_0, z\times a_1, y \times b_0,$ and $y\times b_1$ where $a_i$ are generators of $H^0(O(2))$, $b_i$ generators of $H^1(O(2))$, $z$ a generator of $H^{2n-4}(M_n)$ and $y$ a generator of $H^{2n-5}(M_n).$  Then 
$$Sq^1_t(z\times a_i) = 0 + z\cup w_1 \times a_i = 0$$
and 
$$Sq^1_t(y \times b_i) = z\times b_i + y \cup w_1 \times b_i = 0$$
since one can see by examining the ideal $R$ above that $y \cup w_1 = z.$  Hence $Sq^1_t:\mathbb{H}^0(M_n \times O(2)) \to \mathbb{H}^1(M_n \times O(2))$ is always zero, so if $x \in \mathbb{H}_1(M_n \times O(2))$, then $x Sq^1_t = 0$.  

Now consider $Ad(E_n) = E_n \times _{O(2)} O(2)$, which has the two components $$C_1 := E_n \times _{O(2)} SO(2) \text{  and  } C_2 := E_n \times _{O(2)} ASO(2),$$ where 
$A = \bigl( \begin{smallmatrix} 0 & 1 \\ 1 & 0 \end{smallmatrix} \bigr).$  Now $E_n = V_{2,n}$, the Stiefel manifold of orthonormal $2$-frames in $\mathbb{R}^n$.  So 
$$C_1 = \tilde{Gr}_{2,n} \times _{\mathbb{Z}_2} S^1$$
where $\tilde{Gr}_{2,n}$ is the Grassmannian of oriented $2$-planes and $\mathbb{Z}_2$ is identified with $\{ A, 1 \}$ and so acts on $S^1$ by complex conjugation.  $\tilde{Gr}_{2,n}$ is orientable and the action of $A$ is orientation preserving if and only if $n$ is even; hence since complex conjugation on $S^1$ is orientation reversing, $C_1$ is orientable if and only if $n$ is odd.  

The second component 
$$C_2 = V_{2,n}/ \sim~$$
where, viewing $V_{2,n}$ as the set of $n \times 2$ matrices with orthonormal columns, $B \sim BA \sim - B \sim -BA.$  The map $B \mapsto -B$ is always orientation preserving; the map $B \mapsto BA$ is orientation preserving if and only if $n$ is odd.  Hence $C_2$ is also orientable if and only if $n$ is odd.  

For $n$ even, $M_n$ is orientable but $Ad(E_n)$ is not.  Hence $Sq^1_t = Sq^1$.  By examining each component we see that $Sq^1: \mathbb{H}^0(Ad(E_n)) \to \mathbb{H}^1(Ad(E_n))$ is surjective.  Thus if $x$ is a nonzero element of $\mathbb{H}_1(Ad(E_n))$, then $xSq^1_t \neq 0.$  

Now notice that for any $n$, the map
$$\mathbb{H}_1(Ad(E_{n+1})) \stackrel{\theta _ {i_{n,n+1}}}{\to} \mathbb{H}_1(Ad(E_n))$$
is an isomorphism, since it is $\mathbb{Z}_2$ Poincar\'e dual to $H^0(Ad(E_{n+1})) \stackrel{i_{n,n+1}^*}{\to} H^0(Ad(E_{n})).$  Since $\theta_{i_{n,n+1}}$ also preserves the twisted $\mathcal{A}$ module structure, if $n$ is even and $x$ is a nonzero element of $\mathbb{H}_1(Ad(E_{n+1}))$, then $xSq^1_t \neq 0$.  In particular, since $\mathbb{H}_1(M_{n+1} \times O(2)) \to \mathbb{H}_n(M_n\times O(2))$ is also always an isomorphism, we see that the two pro-ring spectra $\{Ad(E_n)^{-TM_n} \}$ and $\{ (M_n \times O(2))^{-TM_n} \}$ are not homotopy equivalent.

\section{Homotopy Invariance}\label{sec:hominv}
\begin{prop} \cite{cks}
Let $M_1$ and $M_2$ be two smooth, oriented, connected, closed $d$-dimensional manifolds with a smooth, orientation-preserving homotopy equivalence $f:M_1 \to M_2.$  Then the ring spectra maps $$\tilde{f}: LM_1^{-TM_1} \to f^*(LM_2)^{-TM_1}$$ and $$\theta _f: LM_2^{-TM_2} \to f^*(LM_2)^{-TM_1}$$ are both weak equivalences.  Furthermore, $Lf:LM_1 \to LM_2$ is a ring isomorphism in homology.
\end{prop}
 
\begin{proof}
Notice that we have two morphisms of fiberwise monoids
\begin{equation}\label{dia:ftil} \begin{CD}
\Omega M_1 @>\tilde{f}>> \Omega M_2 \\
@VVV @VVV \\
LM_1 @>\tilde{f}>> f^*LM_2 \\
@VVV @VVV \\
M_1 @= M_1
\end{CD} \end{equation}
and 
\begin{equation}\label{dia:f}\begin{CD}
\Omega M_2 @= \Omega M_2 \\
@VVV @VVV \\
f^*LM_2 @>f>> LM_2 \\
@VVV @VVV \\
M_1 @>f>> M_2
\end{CD} \end{equation}
where $\tilde{f}: \gamma \mapsto (f \circ \gamma, \gamma (0)).$
First consider $\tilde{f} : LM_1^{-TM_1} \to f^*LM_2^{-TM_1}.$  Applying the Thom isomorphism $u_*$ in homology to both sides, we obtain the commutative diagram
$$\begin{CD}
H_*(LM_1 ^{-TM_1}) @>{\tilde{f}_*}>> H_*(f^*LM_2^{-TM_1}) \\
@V{u_*}V{\cong}V @V{\cong}V{u_*}V \\
H_{*+d}(LM_1) @>{\tilde{f}_*}>> H_{*+d}(f^*LM_2)
\end{CD}$$
Both $Lf:LM_1 \to LM_2$ and $f: f^*LM_2 \to LM_2$ are homotopy equivalences since $f$ is.  Furthermore $LF = f \circ \tilde{f}$, so $\tilde{f}$ is also a homotopy equivalence. Hence $\tilde{f}_* : H_*(LM_1 ^{-TM_1}) \to H_*(f^*LM_2^{-TM_1})$ is an isomorphism for each value of $*$, so $\tilde{f}$ is a weak equivalence of ring spectra.  

Diagram (\ref{dia:f}) induces the map $\theta _f: LM_2^{-TM_2} \to f^*LM_2^{-TM_1}$, which is the composition 
$$LM_2^{-TM_2}\wedge S \stackrel{id \wedge j}{\to} LM_2^{-TM_2}\wedge {S^N}^{-TS^N} \stackrel{\tau}{\to} f^*LM_2^{-TM_1}.$$
First consider $\tau$.  Applying the cohomology Thom isomorphism $u^*$ on both sides we have the commutative diagram
$$\begin{CD}
H^*(f^*LM_2^{-TM_1}) @>{\tau ^*}>> H^*( LM_2^{-TM_2}\wedge {S^N}^{-TS^N}) \\
@A{u^*}A{\cong}A @A{\cong}A{u^*}A \\
H^{*+d+N}(f^*LM_2^{\nu}) @>{\tau ^*}>> H^{*+d+N}( LM_2\times S^N).
\end{CD}$$
Here the bottom line is the Pontrjagin-Thom collapse map associated to $f\times e: f^*LM_2 \into LM_2 \times S^N,$ where $e:M_1 \into S^N$, and $N>d.$  
Now we can compute $\tau ^*$ via the commutative diagram
$$\begin{CD}
H^*(LM_2 \times S^N) @= H^*(LM_2 \times S^N) \\
@V{\cup t'}VV @VV{u^* \circ (f\times e)^*}V \\
H^{*+N}(LM_2\times S^N, LM_2 \times S^N - (ev\times id)^{-1}(\nu)) @>{\cong}>{excision}> H^{*+N}(f^*LM_2^{\nu}) \\
@V{J}VV @VV{\tau ^*}V \\
H^{*+N}(LM_2 \times S^N) @= H^{*+N}(LM_2 \times S^N)
\end{CD}$$
where $t'$ is the image of the Thom class under the excision isomorphism, $\nu$ is the normal bundle of $M_1$ in $M_2 \times S^N$, and $J$ comes from the relative cohomology long exact sequence.  The left column is a morphism of $H^*(LM_2 \times S^N)$ modules, hence is the cup product with $J(1\cup t')$ where $1\in H^0(LM_2\times S^N)$ (notice that even is $M_2$ is not simply connected, and hence $LM_2$ is not connected, there is still an element $1 \in H^0(LM_2 \times S^N)$ which is an identity for the cup product).  To calculate $J(1\cup t')$ consider the following commutative diagram:
$$\begin{CD} 
H^*(LM_2 \times S^N) @>{u^*\circ (f\times e)^*}>> H^{*+N}(f^*LM_2^\nu) @>{\tau ^*}>> H^{*+N}(LM_2 \times S^N) \\
@A{(ev\times id)^*}AA @AA{ev^*}A @AA{(ev\times id)^*}A \\
H^*(M_2 \times S^N) @>{u^* \circ (f\times e)^*}>> H^{*+N}(M_1^{\nu}) @>{\tau ^*}>> H^{*+N}(M_2 \times S^N).
\end{CD}$$
Now $(ev \times id)^*(1) = 1$ so it's enough to find the image of $1\in H^0(M_2 \times S^N)$ under the bottom row, i.e. $(f\times e)^! \circ (f\times e)^* (1) = (f\times e)^!(1).$  By Poincar\'e duality this is the same as finding 
$$ (f\times e)_* ([M_1]) \in H_d(M_2 \times S^N) \cong H_d(M_2)\otimes H_0(S^N) \cong \mathbb{Z}.$$  Since $f = \pi_{M_2} \circ (f \times e)$ is a homotopy equivalence, $(f\times e)_*([M_1])$ is a unit in $H_d(M_2\times S^N),$ so in particular $(f \times e)^! \circ (f\times e)^* (1) = \pm 1 \otimes z  \in H^0(M_2)\otimes H^N(S^N) \cong H^N(M_2 \times S^N),$ where $z$ is a generator for $H^N(S^N).$  
Hence by naturality of the cross product, 
$$J(1 \cup t') = (ev\times id)^* (f \times e)^!(1) = \pm 1 \otimes z \in 
 H^0(LM_2) \otimes H^N(S^N) \subseteq H^N(LM_2 \times S^N).$$
Hence we see that the composition 
$$\tau ^* \circ u^* \circ (f\times e)^*:H^*(LM_2 \times S^N) \to H^{*+N}(LM_2 \times S^N)$$ 
carries $H^*(LM_2)\otimes H^0(S^N)$ isomorphically onto $H^*(LM_2)\otimes H^N(S^N)$ and is zero on $H^{*-N}(LM_2)\otimes H^N(S^N).$  Furthermore since $f: f^*LM_2 \to LM_2$ is a homotopy equivalence and $e$ is nullhomotopic, $(f\times e)^*$ carries $H^*(LM_2)\otimes H^0(S^N)$ isomorphically onto $H^*(f^*LM_2).$  Thus $\tau^*$ carries $H^{*+N}(f^*LM_2^{\nu})$ isomorphically onto $H^*(LM_2)\otimes H^N(S^N).$

We are now interested in the image of $H^{*+d}(LM_2)\otimes H^N(S^N)$ under 
$$\begin{CD}
H^*(LM_2^{-TM_2}\wedge {S^N}^{-TS^N}) @>{(id\wedge j)^*}>> H^*(LM_2^{-TM_2}\wedge S) \\
@A{u^* \circ \times}AA @A{\cong}A{u^* \circ \times}A \\
H^{*+d}(LM_2)\otimes H^{N}(S^N) @>{id \otimes j'}>> H^{*+d}(LM_2)\otimes H^0(pt.)
\end{CD}$$
where the map $j'$ comes from the unit $j$ of ${S^N}^{-TS^N}$ and the Thom isomorphism:
$$\begin{CD}
H^*({S^N}^{-TS^N}) @>{j^*}>> H^*(S) \\
@A{\cong}AA @AA{\cong}A \\
H^{*+N}(S^N) @>{j'}>> H^*(pt.).  
\end{CD}$$

Now $j$ is Spanier-Whitehead dual to $S^N_+ \to S^0$ and is an isomorphism in cohomology for $*=0$ and is identically zero for $*\neq 0.$  Hence $j'$ is an isomorphism from $H^N(S^N)$ to $H^0(pt)$ and is obviously zero in all other degrees.  In particular, $id \otimes j'$ is an isomorphism from $H^{*+d}(LM_2)\otimes H^N(S^N)$ to $H^{*+d}(LM_2).$  Hence we see that 
$$\theta_f:LM_2^{-TM_2} \to f^*LM_2^{-TM_1}$$ is an isomorphism in cohomology, hence is a weak equivalence of ring spectra.  

We now have weak equivalences 
$$LM_1^{-TM_1} \stackrel{\tilde{f}}{\to} f^*LM_2^{-TM_1}$$
and 
$$LM_2^{-TM_2} \stackrel{\theta _f}{\to} f^*LM_2^{-TM_1}.$$
Applying the Thom isomorphism on both sides, we obtain coalgebra isomorphisms
$$H^*(f^*LM_2) \stackrel{\tilde{f}^*}{\to} H^*(LM_1)$$
and
$$H^*(f^*LM_2) \stackrel{\theta _f^*}{\to} H^*(LM_2).$$
Now the computations above show that the composition
$$H^*(LM_2)\stackrel{\pi _{LM_2}^*}{\to} H^*(LM_2\times S^N) \stackrel{(f\times e)^*}{\to} H^*(f^*LM_2) \stackrel{\theta _f^*}{\to}  H^*(LM_2)$$ 
is the identity (since it is the cup product with $1$).  
Since $\pi_{LM_2} \circ (f\times e) = f : f^*LM_2 \to LM_2$, this shows that $\theta _f^* \circ f^*$ is a coalgebra isomorphism.  Since $\theta _f^*$ is also, $f^*$ must be a coalgebra isomorphism. Finally, $Lf = f \circ \tilde{f} :LM_1 \to LM_2$, so $(Lf)^*$ is a coalgebra isomorphism in cohomology.  It follows from the universal coefficient theorem that $f_* \circ {\theta _f}_*$ is the identity on $H_*(LM_2)$; hence $f_*$ and therefore $Lf_*$ are ring isomorphisms.  

\end{proof}

\section{Examples and applications} \label{examples}

{\em Mapping spaces}

\begin{prop} \cite{SK}

Let $M$ be a closed manifold and $N$ a closed $n$-manifold. Then

\

a) $map(S^n,M)^{-TM}$ is an $E_n$-ring spectrum;

\

b)  
$map(N,M)^{-TM}$ is a module
over  $map(S^n,M)^{-TM}$    .

\end{prop}

\begin{proof} We must show that $map(S^n,M) \to M$ is a fiberwise $C_n$-space,
and 

\noindent $map(N,M) \to M$ a module over it.
Let us consider the action on fibers, that is classical.
A little $n$-cube in $\mathcal{C}_n(i)$
 is an affine embedding $\coprod_i I^n \to I^n$ of a union of
cubes into a cube. The one-point compactification $S^n \to \vee_i S^n$ defines
by precomposition the classical action $\CC_n(i) \x (\Omega^n M)^i \to \Omega^n M$ on
a based $n$-fold loop space.

For the module action, take
a closed $n$-disc $\bar{D}$ embedded into $N$, with interior $D$.
We get a collapse map $N \to N/\bar{D} \vee N/(N-D)$. If we compose by
standard identifications $N/\bar{D} \cong N$ and $N/(N-D)=\bar{D}/(\bar{D}-D)\cong S^n$,
then we recover the coaction $N \to N \vee S^n$ associated to the cofibration
sequence of the top cell attaching map $S^{n-1} \to N_{n-1} \to N$.
If we precompose by the coaction we obtain the action
$\Om^n M \x map_*(N,M) \to map_*(N,M)$. This gives a 
module structure of $map_*(N,M)$ over the $\mathcal{C}_n$-algebra
$\Om^n M$ (up to coherent homotopies) in the strong sense of
\cite{gk}.
The fiberwise structure is defined by precomposing 
with the coaction map, as in the based case,
 except that
the base point in $M$ is now allowed to vary, according to the fiber \cite{SK}.

\end{proof}

In \cite{SK} the string product is used
to compute the homology of $map(S^2,\CP^n)$.
Moreover the module structure gives a periodicity result 
about the homology of $map(\Sigma,\CP^n)$ when $\Sigma$ is a Riemann surface.
The homology mod $p$ depends only on the residue mod $p$ of the component degree
in $\Z=\pi_0(map(\Sigma,\CP^n))$.

\

 {\it Universal example}:

 Let $Y$ be a topological monoid (or a $C_n$-space) with an action of a topological
 group $G$ preserving the sum,
 and $G \to E \to M$ a $G$-principal fibration 
  over a closed manifold $M$.
Then $Y \to Y\times_G E \to M$ is a fiberwise monoid (or a $C_n$-space) .
This will be useful in the next few examples.

The example is universal because any fiberwise monoid has this form with
$G=\Omega M$ the based loop space and $E=PM$ the contractible path space.
Up to homotopy we can always replace $\Omega M$ by a topological group.

\

{\em Rational curves}

 Let $M$ denote a homogeneous complex projective manifold. Then $M=G/P$,
where $G$ is a complex semisimple Lie group and $P$ a parabolic subgroup.
Let $hol(G/P)$ and $hol_*(G/P)$ be respectively the spaces of free
 and based holomorphic maps from the Riemann sphere to $G/P$.
\begin{prop}   \cite{SK}
The Thom  spectrum $hol(M)^{-TM}$ is an $E_2$-ring spectrum.
\end{prop}
\begin{proof}
It is sufficient to show that the bundle $hol(G/P) \to G/P$ has a fiberwise
 action of the little 2-cubes operad.

In \cite{BHMM} $hol_*(G/P)$ is identified to a suitable
 configuration space of particles in the plane (or a disc) with labels.
 The little 2-cubes operad $\mathcal{C}_2$ acts on $hol_*(G/P)$
by its standard action on the locations of the particles
 and keeping the labels unchanged.
In the updated version of \cite{SK} 
we explain that the operad action is $U$-equivariant, with $U \subset P$
and $V \subset G$ suitable maximal compact subgroups such that $U \subset V$.
Then the space
$hol(G/P)=V \x_U hol_*(G/P)$, according to the universal example,
has a fiberwise action of the
little 2-cubes operad.

\end{proof}

In \cite{SK} the homology of $hol(\CP^n)$ was computed using
the string product.

\

{\em Remark}: There is no module structure in the holomorphic category as opposite to the continuous:
the space of holomorphic based maps
$hol_*(\Sigma,G/P)$ from a Riemann surface $\Sigma$ to $G/P$
is not a module over $hol_*(G/P)$.

\

{\em Spaces of knots and embeddings}

\begin{prop}
Let $K$ be the space of long knots in $\R^3$ and $Emb(S^1,S^3)$ the
space of knots in $S^3$.
Then there is a map of $E_2$-ring spectra
$\Sigma^{\infty-5} Emb(S^1,S^3)_+ \to \Sigma^{\infty} K_+$.
\end{prop}

A long knot is an embedding $f:\R^1 \to \R^3$ such that
$f(t)=(t,0,0)$ for $|t|\geq 1$.
The operad action on the fiber $K$ was discovered by
Budney \cite{Budney}.
He constructs a space $EC(1,D^2)$ homotopy equivalent to $K \times \Z$.
This is the space of embeddings $f:\R \times D^2 \to \R \times D^2$
with support contained in $I \x D_2$, where $D^2$ is the standard
unit 2-disc. The operad action is then constructed by rescaling
the embeddings in the $\R$-variable and composing them in a suitable
order. The rescaling depends on the first coordinate of the little 2-cubes,
and the order depends on the second coordinate.
The group $SO(2)$ acts on $EC(1,D^2)$ by conjugation, namely
$\theta(f)=(id_\R \x \theta) \circ f \circ (id_\R \x \theta^{-1})$.
It is easy to see that the operad composition is $SO(2)$-equivariant
with respect to this action. This implies that
$SO(4) \x_{SO(2)} EC(1,D^2) \to SO(4)/SO(2)$ has a fiberwise action of
$\mathcal{C}_2$. Let $EC_0(1,D^2)$ be the component of the identity in $EC(1,D^2)$.

\

The forgetful map
$EC_0(1,D^2) \to K$ sending $f$ to the long knot $t \mapsto  f(t,0,0)$
is a homotopy equivalence \cite{Budney} but also a $SO(2)$-equivariant
map. The $SO(2)$ action on $K$ is induced by the standard action on the two last coordinates
of $\R^3$. This gives a homotopy equivalence
$SO(4) \x_{SO(2)} EC_0(1,D^2) \simeq SO(4) \x_{SO(2)} K$. There is also a
homotopy equivalence $SO(4) \x_{SO(2)} K \simeq Emb(S^1,S^3)$
\cite{bu-co}. This gives a $E_2$-ring spectrum structure on $Emb(S^1,S^3)^{-T(SO(4)/SO(2))}$,
but $SO(4)/SO(2) \cong S^3 \x S^2$ has a trivial normal bundle. The map in the statement is induced by
the inclusion $K \to Emb(S^1,S^3)$ by Proposition \ref{en-maps}.

\

{\em Gauge groups}

Let $G$ be a compact Lie group.
Let $\G(N)$ be the topological
groupoid with one object for any isomorphism class of  $G$-principal bundles over $N$,
and a morphism for each $G$-equivariant map covering the identity of $N$.
The classifying space $B\G(N)$ is the union of the classifying spaces of
all gauge groups of a bundle over all isomorphism types of bundles.

\begin{prop}
For any $n$ there is a pro-$E_n$-ring spectrum $B\G(S^n)^{-TBG}$ .
For any closed $n$-manifold $N$ there is a pro-module spectrum
$B\G(N)^{-TBG}$ over 

\noindent $B\G(S^n)^{-TBG}$.
\end{prop}

The construction for $n=1$ recovers the string topology of $BG$ defined
in Section \ref{sec:bg}.  
We recall first that there is a homotopy equivalence
$B\G(N) \simeq map(N,BG)$ \cite{gauge}.
We consider the fibration $p:map(N,BG) \to BG$, a filtration
by closed submanifolds $BG=colim_{m \to \infty} (M_m)$, and define
$X_m(N)=p^{-1}(M_m)$. Then, for any $m$, $X_m(S^n)$ is a fiberwise monoid over $M_m$
with fiber $\Omega^n(BG)$. Consequently
$X_m(S^n)^{-TM_m}$ is a $E_n$-ring spectrum.
 Moreover $X_m(N)$ is a fiberwise module
over $M_m$ with fiber $map_*(N,BG)$, and thus $X_m(N)^{-TM_m}$ is
a module over $X_m(S^n)^{-TM_m}$.
We conclude by defining as in Section \ref{sec:bg}
the pro-$E_n$-ring spectrum $B\G(S^n)^{-TBG}$ 
$$ \dots \to   X_{m+1}(S^n)^{-TM_{m+1}}\to  X_m(S^n)^{-TM_m} \to \dots \to X_1(S^n)^{-TM_1}$$
and its pro-module spectrum  $B\G(N)^{-TBG}$ 
$$ \dots \to   X_{m+1}(N)^{-TM_{m+1}}\to  X_m(N)^{-TM_m} \to \dots \to X_1(N)^{-TM_1}.$$

\

In \cite{SS} the homology operations will be used to determine
the homology of $B\G(S^4)$ for $G=SU(2)$.
 This technique shows that there is no periodicity
 pattern, and the homology mod $p$ of $B\G_k(N)$, the component corresponding to 
 bundles with second Chern number $k$, gets more complicated
as the $p$-exponent of $k$ grows.

\

{\em Maps from Riemann surfaces}

\begin{prop}
For any non-negative integer $g$ let $F_g$ be a Riemann surface of genus $g$.
Let $M$ be a closed manifold. Then the wedge
$\vee_{g \in \N}map(F_g,M)^{-TM}$ is a ring spectrum.
\end{prop}

We show that the union $\coprod_{g \in \N}  map(F_g,M)$ is a fiberwise monoid over
$M$.

For any genus $g$ we have a cofibration sequence
$S^1 \to \vee_{2g} S^1 \to F_g$. The first map represents in
the fundamental group $\pi_1(\vee_{2g}S^1)$, free on generators
$y_1,\dots,y_{2g}\,$, the product of commutators $[y_1,y_2] \dots [y_{2g-1},y_{2g}]$.

We may suppose that $\Omega M$ is a group by taking
for example the realization of the Kan loop group on the singular
simplicial set of $M$. 
Thus we can identify $map_*(F_g,M)$
to the homotopy fiber of $\alpha_g: \Omega M^{2g} \to \Omega M$,             
where $\alpha_g(x_1,\dots,x_{2g}) = \prod_{i=1}^g [x_{2i-1},x_{2i}]$.

Now the diagram 
$$\begin{CD} \label{comm}
(\Om M)^{2g} \x (\Om M)^{2h} @= (\Om M)^{2(g+h)} \\ 
@VV{\alpha_g \x \alpha_h}V       @VV{\alpha_{g+h}}V  \\     
\Om M        \x  \Om M     @>m>>  \Om M
\end{CD} $$

strictly commutes, where $m$ is the loop sum.
On homotopy fibers this induces maps
$$\mu_{g,h} : map_*(F_g,M) \times map_*(F_h,M) \to map_*(F_{g+h},M)$$ making
the disjoint union $\coprod_g map_*(F_g,M)$ into a monoid.
Now $\Omega M$ acts on itself by conjugation. With respect to this
action on each coordinate all maps in the diagram above
 are $\Omega M$-equivariant.
Thus $map_*(F_g,M)$ inherits an action of $\Omega M$ and the multiplication
maps $\mu_{g,h}$ are  equivariant. 
As a consequence the union  over all $g's$ of the
homotopy orbit spaces 
$PM \times_{\Om M} map_*(F_g,M)$
 is a fiberwise monoid over $PM \times_{\Om M} * \simeq M$.
But the action of $\Om M$ is the
usual action moving the base point around
loops, and it is well known that with respect to this action
 $PM \x_{\Om M} map_*(F_g,M) \simeq
map(F_g,M)$.

\end{document}